\documentclass[12pt]{article}
\usepackage{latexsym}
\usepackage{amssymb}

\begin{document}
\begin{center}


{\bf \Large Domains of holomorphy of generating functions of
P\'olya frequency sequences of finite order.   }
\end{center}

\begin{center}
Maria Teresa Alzugaray (mtrodrig@ualg.pt)

\end{center}

{\bf Abstract:} A domain $G\subset \overline\mathbb C$ is the
domain of holomorphy of the generating function of a P\'olya
frequency sequence of order $r$ if and only if it satisfies the
following conditions: (A) $G$ contains the point $z=0$, (B) $G$ is
symmetric with respect to the real axis, (C)  $T=dist(0,\partial
G)\in
\partial G$.

\section{Introduction.}

The P\'olya frequency sequences, also called multiply positive sequences,
were first introduced by Fekete
in 1912 (see~\cite{fek}). They were studied in detail by Karlin
(see~\cite{tp}).

The class of all P\'olya frequency  sequences of order
$r\in \mathbb{N}\cup \{\infty\}$ ($r$-times positive) is denoted by
$PF_r$ and consists of the sequences $\{c_k\}_{k=0}^\infty$ such that all
minors of order $\le r$ (all minors if $r=\infty$) of the infinite matrix

 $$
 \left\|
  \begin{array}{ccccc}
   c_0 & c_1 & c_2 & c_3 &\ldots \\
   0   & c_0 & c_1 & c_2 &\ldots \\
   0   &  0  & c_0 & c_1 &\ldots \\
   0   &  0  &  0  & c_0 &\ldots \\
   \vdots&\vdots&\vdots&\vdots&\ddots
  \end{array}
 \right\|
 $$
are non-negative. The class of corresponding generating functions
$$
f(z)=\sum_{k=0}^\infty c_kz^k
$$
is also denoted by $PF_r$. The radius of convergence of a $PF_r$ generating
function ($PF_r$ g.f.) is positive
provided $r\ge2$ (\cite{tp}, p.394). Further we will suppose,
without loss of generality,
that $c_0=1$.

The class $PF_\infty$ was completely described in~\cite{aissen}
(see also \cite{tp}, p. 412):

{\bf Theorem~\cite{aissen}:} {\it The class $PF_\infty$ is formed by the
functions
$$
f(z)=e^{\gamma z}\prod_{k=1}^\infty (1+\alpha_kz)/(1-\beta_kz),
$$
where $\gamma\ge0,\alpha_k\ge0,\beta_k\ge0,\sum(\alpha_k+\beta_k)<\infty.$}
\medskip

In 1955, I.J. Schoenberg set up the problem of characterizing the classes
$PF_r, r\in \mathbb{N}$. Some results were obtained
that showed essential differences between the properties
of $PF_\infty$ g.f.
and those  of $PF_r$ g.f. with $r\in \mathbb{N}$
(see~\cite{katos,ind} and~\cite{essing}).

\section{Statement of results.}

This paper deals with the study of  $PF_r$ domains of holomorphy
with $ r\in \mathbb{N}$.

A domain $G\subset\overline{\mathbb{C}}$
is called a {\it $PF_r$ domain of holomorphy} if there exists a
$PF_r$ g.f. which is analytic in $G$ and admits no analytic
continuation across the boundary of $G$.

It follows from Theorem~\cite{aissen} that if $G$ is a $PF_\infty$
domain of holomorphy,   then ${\mathbb{C}}\backslash G$ is at most a countable set of
points $\{p_k\}$ on the positive ray  such that $\sum (1/p_k)<\infty$
(the points cannot be "too close" to each other).

The situation with the $PF_r$ domains of holomorphy for which
$r\in \mathbb{N}$,
is quite different.
They can be much more complicated as the main result of~\cite{sing} shows:

{\bf Theorem~\cite{sing}:} {\it Let $E$ be a closed set in $\mathbb{C}$,
satisfying the
conditions: (i) $E$ is symmetric with respect to the real axis,
(ii) $E \cap\{z:|z|\le 1\}=\emptyset$. For any $r\in \mathbb{N}$,
there exists a
function $f(z)$ such that: (i) $f(z)\in PF_r$, (ii) the set of all
singularities of $f(z)$ coincides with $E\cup\{1\}$.}

From the theorem above, we obtain conditions,
which, taken together, are sufficient
for a domain $G$
to be a  $PF_r$ domain of holomorphy with $r\in \mathbb{N}$:

(I) $G$ contains the point $z=0$;

(II) $G$ is symmetric with respect to the real axis;

(III) Let $E$ be the set from Theorem 3. $\partial G=\{T\}\cup
\partial E$ with $T\in{\mathbb{R}}$, $T>0$ and $E \cap\{z:|z|\le
T\}=\emptyset$.

On the other hand, it is not difficult to see that if $G$ is
the domain of holomorphy of the $PF_r$ g.f.
$f(z)=\sum _{k=0}^{\infty}c_kz^k, 2\le r<\infty,$ then the following
conditions
must be satisfied:

(A) $G$ contains the point $z=0$  (this condition is assured by
$f(z)\in PF_r\subset PF_2$);

(B) $G$ is symmetric with respect to the real axis (since
$c_k\in{\mathbb{R}}, k=0,1,2,...$);

(C)  $T=dist(0,\partial G)\in \partial G$
(by the well-known Pringsheim's theorem
on singularities of power series with non-negative coefficients).

The class of domains satisfying (A)-(C) is much larger than the
class satisfying (I)-(III).
The question of the exact description
of all $PF_r$ domains of holomorphy arises.

The following result shows that condition (III) is not necessary, since
$\partial G$ can contain other points of $\{z:|z|=T\}$ besides $T$.

{\bf Theorem 1.} {\it Let $g(z)=\sum_{k=0}^\infty b_kz^k$ be a
function with real and bounded Taylor coefficients, i.e.

\begin{equation}
\label{lim}
|b_k|<C,\  k=0,1,2,\ldots,
\end{equation}
for some constant $C=C(g),C>0.$ Then for any $r\in \mathbb{N}$
there exists a $\varepsilon >0$
such that the function

$$
f_{\varepsilon}(z)=\frac{1}{(1-z)^{r^2}}+\varepsilon g(z)
$$
is a  $PF_r$ g.f..}

The function $h(z)=\sum_{k=0}^\infty z^{k!}$ satisfies (\ref{lim})
so for  $g=h$ the set $\{z:|z|=1\}$ coincides with the
singularity set of $f_{\varepsilon}$, i.e. $\mathbb{D}=\{z:|z|<1\}$ is a $PF_r$
domain of holomorphy (from now on we denote the unit disc by $\mathbb{D}$).

{\bf Theorem 2:} {\it The domain $\Omega$ is the domain of
holomorphy of a function with bounded Taylor coefficients if and
only if $ \Omega$ contains the unit disc.}

The necessity of the condition $\mathbb{D}\subset\Omega$ is evident.

With the help of Theorem 2 we arrive at the main result of our paper:

{\bf Theorem 3:} {\it The domain $G$ is a $PF_r$ domain of
holomorphy if and only if $G$ satisfies conditions} (A)-(C).

It follows from our theorems  and the proof of Theorem 2 that
among the $PF_r, r\in\mathbb{N},$ there are functions that have no limit
near the points of the boundary of their domains.
Corollary 1 of the following theorem shows that a $PF_r$ g.f.
with $r\ge 2$ cannot be bounded in
its domain of holomorphy.

{\bf Theorem 4:} Let $f(z)$ be a $ PF_r$ g.f., $r\ge 2$,
and $T$ be its radius of convergence, $T<\infty$. Then
$ (1-z/T)f(z)$ is a $PF_{r-1}$ g.f..

{\bf Corollary 1:} Let $f(z)$ satisfy the conditions of Theorem 4
for $r=2$. Then
$$
\lim_{x\to T^-}(1-x/T) f(x)=a>0,
$$
where $a$ can be equal to $+\infty$.

{\bf Corollary 2:} Let $f(z)$ satisfy the conditions of Theorem 4
and $T$ be an essential singularity of the function $f$. Then
$$
\lim_{x\to T^-}(1-x/T)^{r-1} f(x)=a>0,
$$
where $a$ can be equal to $+\infty$.

\section{ Proofs of the results:}

For proving Theorem 1 we will need two lemmas whose proofs
can be found in \cite{sing}:

{\bf Lemma 1:} {\it Let $r,n \in {\mathbb{N}}, r \ge 1$ and $1 \le n \le r.$
 Let

 \begin{equation}
 \frac{1}{(1-z)^r}=\sum_{k=0}^\infty
 \frac{(k+1)\cdots(k+r-1)}{(r-1)!}z^k
 =\sum_{k=0}^\infty a_k(r)z^k  \label{coef}
 \end{equation} and

 \begin{equation}a_k(r)=0\ \ { for}\ \ k<0\label{coef0}.\end{equation}
 Then for $k\ge 0,$
 $$
 \det\| a_{k+j-i}\| _{i,j=\overline{1,n}}=
 \prod_{i=1}^n \frac{(i-1)!}{(r-i)!} (k+i)\cdots (k+i+r-n-1),
 $$
 where we consider $(k+i)\cdots (k+i+r-n-1)=1$ for $n=r.$}

{\bf Lemma 2:} (see \cite{schoenb}) {\it
If $ \sum_{k=0}^\infty c_k  \lambda ^k<\infty $ for certain $\lambda$,
 $ 0<\lambda <1$, $c_k=0$ for $k<0$
 and $ \det\| c_{k+j-i}\| _{i,j=\overline{1,n}}>0 $
 for all $k \ge 0$ and $n, 1\le n\le r$, then
 $\{ c_k \}_{k=0}^\infty \in PF_r.$}

Now, we are going to prove the following fact,
which is slightly more general
than Theorem 1.

\medskip

{\bf Theorem 1':} Let the function $g$ satisfy the hypothesis of
Theorem 1. Then for any $r\in N$ and any $\alpha\in\mathbb{N}$,
there exists a $\varepsilon >0$ such that the functions
$f_{\varepsilon}(z),f'_{\varepsilon}(z),\ldots,
f^{(\alpha)}_{\varepsilon}(z)$ are $PF_r$ g.f..

\medskip
{\bf Proof of Theorem 1':} Let $r\in\mathbb N$ be a fixed number.
Let $f^{(p)}_\varepsilon (z)=\sum _{k=0}^\infty c^p_k(r^2) z^k$,
where $0\le p\le \alpha$. Then $c^p_k(r^2)=a^p_k(r^2)+\varepsilon
b^p_k$ where $a^p_k(r^2) $ and $b^p_k$ are the Taylor coefficients
of the functions

$$
\frac{r^2\cdots(r^2+p-1)}{(1-z)^{r^2+p}}
\ \ \left(\frac 1{(1-z)^{r^2}}\ \ \mbox{for}\ \ p=0\right)
$$
and $g^{(p)}(z)$ respectively (consider $a_k^0(r^2)=a_k(r^2)$ and
$b_k^0=b_k$). More explicitly,

\begin{equation}\label{der}
\begin{array}{ll}
a^p_k(r^2)=\frac{(k+1)\cdots (k+r^2+p-1)}{(r^2-1)!},&{\rm for}\
r-1+p>0,k=0,1,2,\ldots;\\
a_k^0(1)=1,& {\rm for}\ k=0,1,2,\ldots;\\
a_k^p(r^2)=0,& {\rm for}\ k<0; \\
\end{array}\end{equation}

and
\begin{equation}
\begin{array}{ll}
b^p_k=(k+p)\cdots(k+1)b_{k+p},& {\rm for}\ k=0,1,2,\ldots,p\neq 0;\\
b_k^0=b_k,& {\rm for}\ k=0,1,2,\ldots;\\
b^p_k=0,& {\rm for}\ k<0.\end{array} \label{bound}\end{equation}
Therefore, for each $p, 0\le p\le\alpha,$

\begin{equation}
\label{ots} a_k^p(r^2)=O(k^{r^2+p-1})\ \ \mbox{and}\ \
b^p_k=O(k^{p}),\ k\to \infty.
\end{equation}

For $n,1\le n\le r$, we have
$$ \left|
\begin{array}{cccc}
c^p_k&c^p_{k+1}&\ldots&c^p_{k+n-1}\\
c^p_{k-1}&c^p_k&\ldots&c^p_{k+n-2}\\
\vdots&\vdots&\ddots&\vdots\\
c^p_{k-n+1}&c^p_{k-n+2}&\ldots&c^p_{k}
\end{array}
\right|=
$$
$$
=
\left|
\begin{array}{cccc}
a^p_k+\varepsilon b^p_k&a^p_{k+1}+\varepsilon b^p_{k+1}&\ldots&a^p_{k+n-1}+\varepsilon b^p_{k+n-1}\\
a^p_{k-1}+\varepsilon b^p_{k-1}&a^p_k+\varepsilon b^p_k&\ldots&a^p_{k+n-2}+\varepsilon b^p_{k+n-2}\\
\vdots&\vdots&\ddots&\vdots\\
a^p_{k-n+1}+\varepsilon b^p_{k-n+1}&a^p_{k-n+2}+\varepsilon b^p_{k-n+2}&\ldots&a^p_k+\varepsilon b^p_k
\end{array}
\right|= $$
$$= \sum_{l=0}^n\varepsilon^lS^p_l(k,n), \ \mbox{where}
\ S^p_0(k,n)=
\left|
\begin{array}{cccc}
a^p_k&a^p_{k+1}&\ldots&a^p_{k+n-1}\\
a^p_{k-1}&a^p_k&\ldots&a^p_{k+n-2}\\
\vdots&\vdots&\ddots&\vdots\\
a^p_{k-n+1}&a^p_{k-n+2}&\ldots&a^p_{k}
\end{array}
\right|,
$$
$$
S_n^p(k,n)=
\left|
\begin{array}{cccc}
b^p_k&b^p_{k+1}&\ldots&b^p_{k+n-1}\\
b^p_{k-1}&b^p_k&\ldots&b^p_{k+n-2}\\
\vdots&\vdots&\ddots&\vdots\\
b^p_{k-n+1}&b^p_{k-n+2}&\ldots&b^p_{k}
\end{array}
\right| \qquad \mbox {and}
$$
$S^p_l(k,n),1 \le l \le n-1 $, is the sum of determinants
of all possible matrices formed by $n-l$ rows from $S^p_0(k,n)$ and $l$
complementary rows from $S^p_n(k,n)$.

The sum $S^p_l(k,n)$ has $n \choose l$ summands. Each summand is a determinant
of order $n$. Each determinant is the sum of $n!$ summands which are the
products of $n-l$ entries of the matrix of $S^p_0(k,n)$ and $l$ entries of
the matrix $S^p_n(k,n)$. It follows from (\ref{ots}) that the modulus of each entry
of $S^p_n(k,n)$ does not exceed $C_p(k+1)^p$, where $C_p$ is a
constant, and each element of
$S^p_0(k,n)$ is non-negative and does not exceed $a^p_{k+r-1}$, since
$\{a^p_k\}$ is a non-decreasing sequence. Hence, we have for $l,1\le l\le n$,
$$
|S^p_l(k,n)|<{n\choose l}n!(C_p)^l(k+1)^{pl}(a^p_{k+r-1})^{n-l} \le
B(k+1)^{pl}(a^p_{k+r-1})^{n-l},
$$
where $B=2^rr!\max_{0\le p\le \alpha}(C_p)^r$. Obviously $B$ does not
depend on $k$.

As a result, we have
$$ \left|
\begin{array}{cccc}
c^p_k&c^p_{k+1}&\ldots&c^p_{k+n-1}\\
c^p_{k-1}&c^p_k&\ldots&c^p_{k+n-2}\\
\vdots&\vdots&\ddots&\vdots\\
c^p_{k-n+1}&c^p_{k-n+2}&\ldots&c^p_{k}
\end{array}
\right| \ge
$$
$$
\ge
S^p_0(k,n)-\sum_{l=1}^n\varepsilon^l|S^p_l(k,n)|\ge
$$
\begin{equation}
\ge S^p_0(k,n)-\varepsilon
B\sum_{l=1}^n(k+1)^{pl}(a^p_{k+r-1})^{n-l}
\label{desig}\end{equation} for $\varepsilon <1$.

From (\ref{ots}) we have
$$
(k+1)^{pl}(a^p_{k+r-1})^{n-l}= O(k^{(r^2+p-1)(n-l)+pl})=
O(k^{(r^2-1)(n-l)+np}),
k\to \infty.
$$

By (\ref{coef}), (\ref{coef0}) and (\ref{der}) we have
$$
a_k^p(r^2)=\frac {r^2+p-1}{(r^2-1)!}a_k(r^2+p),\ k\in\mathbb Z.
$$
Hence, using Lemma 1 we obtain

$$
S^p_0(k,n)= \prod_{j=0}^{p-1}(r^2+j)^n
\prod_{i=1}^n\frac{(i-1)!}{(r^2+p-i)!}(k+i)\cdots(k+i+r^2+p-n-1)
$$
\begin{equation}\label{desig2} \ge M(k+1)^{n(r^2+p-n)},\ n=1,2,\ldots,r,\
k=0,1,2,\ldots,\end{equation} where $M$ is a positive number.

Also, for each $n,1\le n\le r,$ and $l, 1\le l\le n,$
$$
n(r^2+p-n)=n(r^2-1)-(n^2-n)+np\ge n(r^2-1)-(n^2-1)+np
$$
$$
\ge n(r^2-1)-(r^2-1)+np\ge n(r^2-1)-l(r^2-1)+np=(r^2-1)(n-l)+np.
$$

Hence, using (\ref{ots}), for each $p, 0\le p\le\alpha,$ and for
any $n, 1\le n\le r,$

$$
B\sum_{l=1}^n(k+1)^{pl}(a^p_{k+r-1})^{n-l}=O(k^{n(r^2+p-n)}),\ k\to \infty
$$
and by (\ref{desig2}) there is a $\varepsilon^p_n>0$ such that the
inequality
$$
\varepsilon_n^p
B\sum_{l=1}^n(k+1)^{pl}(a^p_{k+r-1})^{n-l}<
\frac12S^p_0(k,n)
$$
holds for any $k \ge 0$.

Let $\varepsilon =\min\{\varepsilon_n^p\}$. Since Lemma 1
$S_0^p(k,n)>0$, then, using (\ref{desig}), for any $k\ge 0$, any
$n,1\le n\le r$, and any $p,0\le p\le\alpha$ we have
$$
\left|
\begin{array}{cccc}
c^p_k&c^p_{k+1}&\ldots&c^p_{k+n-1}\\
c^p_{k-1}&c^p_k&\ldots&c^p_{k+n-2}\\
\vdots&\vdots&\ddots&\vdots\\
c^p_{k-n+1}&c^p_{k-n+2}&\ldots&c^p_{k}
\end{array}
\right|
>\frac12 S^p_0(k,n)>0.
$$

Relying on Lemma 2 we conclude that $\{c^p_k\}_{k=0}^\infty \in
PF_r.$ $\Box$

{\bf Proof of Theorem 2:}

The necessity of the condition $\mathbb{D}\subset\Omega$ is evident.

The sufficiency will be proved by constructing a function $g(z)$
with bounded Taylor coefficients, analytic in $\Omega$,
that cannot be analytically continued through $\partial\Omega$.

 Let $\Omega$ be a domain with $\mathbb{D}\subset \Omega$ and
$\{\zeta_k\}_{k=1}^\infty\subset\partial\Omega$ a
countable set of points dense in $\partial\Omega$.

In the neighborhood of radius $1/2n, n\in \mathbb{N},$ of the point $\zeta_k$
there exists a point $z(n,k)\in \Omega$. Let us denote by
$\lambda(n,k)$ the closest to $z(n,k)$ point of $\partial \Omega$,
i.e. $|\lambda(n,k)-z(n,k)|=dist(z(n,k),\partial \Omega)$.
It is obvious that $|\zeta_k-\lambda(n,k)|<1/n, n\in \mathbb{N}$.
Hence, the countable set $\{\lambda(n,k)\}_{n,k=1}^\infty$
is dense in $\partial\Omega$.
Let us number the sets $\{z(n,k)\}$ and $\{\lambda(n,k)\}$ with
one parameter $k\in \mathbb{N}$ preserving the correspondence
$|z_k-\lambda_k|=dist(z_k,\partial\Omega)$.
Now we choose a sequence $\{d_k\}_{k=1}^\infty\subset \mathbb{R}_+$
such that $\sum_{k=1}^\infty d_k<\infty$ and put

$$
g(z)=\sum_{k=1}^\infty \frac{d_k}{\lambda_k-z}.
$$

If $z\in K$, $K$ is a compact  in $\Omega$, and
$\delta=dist(K,\partial\Omega)$, then
$$
\left|\frac {d_k}{\lambda_k-z}\right|<\frac{d_k}{\delta}.
$$
Thus, the series defining $g$ converges uniformly on
each compact $K\subset\Omega$ and, therefore,
$g$ is analytic in $\Omega$.

Now, we prove that $g$ cannot be analytically continued through
$\partial\Omega$ applying an idea that can be found in Levin's
\cite{lev}, page 117. We fix a point $\lambda_p$ and will show
that $g(z)$ tends to infinity for certain $z$ approaching
$\lambda_p$. Let $N_p$ be a number such that
$$
\sum_{k=N_p+1}^\infty d_k<\frac{d_p}2.
$$
Then
$$
|g(z)|\ge \frac{d_p}{|\lambda_p-z|}-
\left|\sum_{k=1,k\ne p}^{N_p}\frac{d_k}{\lambda_k-z}\right|
-\sum_{k=N_p+1}^{\infty}\frac{d_k}{|\lambda_k-z|}.
$$

Note that for $z=\alpha z_p+(1-\alpha)\lambda_p,0<\alpha<1,$ and
$k\ne p$ the inequalities
$$
|\lambda_k-z|>|\lambda_k-z_p|-|z-z_p|
\ge |\lambda_p-z_p|-|z-z_p|=|\lambda_p-z|
$$
hold.

Hence,
$$
|g(z)|\ge \frac{d_p}{2|\lambda_p-z|}-
\left|\sum_{k=1,k\ne p}^{N_p}\frac{d_k}{\lambda_k-z}\right|,
$$
which shows that $g(z)\to\infty$ when $z\to\lambda_p,z=\alpha
z_p+(1-\alpha)\lambda_p,0<\alpha<1$.

We still have to show that $g$ is a function with bounded Taylor
coefficients.
Denoting
$$
g(z)=\sum_{n=0}^\infty b_nz^n
$$
we have
$$
b_n=\frac{g^{(n)}(0)}{n!}=
\sum_{k=1}^\infty\frac{d_k}{\lambda_k^{n+1}},\ \ {\rm for}\ \  n=0,1,2,\ldots .
$$

But, $|\lambda_k|\ge 1$ and $\sum_{k=1}^\infty d_k<\infty$,
so the sequence $\{b_n\}_{n=0}^\infty$ is bounded.
$\Box$

Theorem 3 follows at once from Theorem 1 and Theorem 2.

{\bf Proof of Theorem 4:}
Since a $PF_2$ g.f. is either a polynomial  or a trascendental function
with positive Taylor coefficients (see~\cite{tp}, page 393), we only
consider the case $c_k>0$ for $k=0,1,2,\ldots$, and $c_k=0$ for $k<0$.

Following~\cite{aissen} (see also~\cite{tp}, page 407), we put
$$
f_1(z)=(1-z/T)f(z)=\sum_{k=0}^\infty c_k^{(1)}z^k,
$$
where $c_k^{(1)}=c_k-c_{k-1}/T$.

Every minor
$$
\left|
\begin{array}{ccccc}
c_{k_1}&c_{k_2}&\ldots&c_{k_n}&\frac{c_{m}}{c_{m-n}}\\
c_{k_1-1}&c_{k_2-1}&\ldots&c_{k_n-1}&\frac{c_{m-1}}{c_{m-n}}\\
\vdots&\vdots&\ &\vdots&\vdots\\
c_{k_1-n}&c_{k_2-n}&\ldots&c_{k_n-n}&\frac{c_{m-n}}{c_{m-n}}\\

\end{array}
\right|
=\frac1{c_{m-n}}
\left|
\begin{array}{ccccc}
c_{k_1}&c_{k_2}&\ldots&c_{k_n}&{c_{m}}\\
c_{k_1-1}&c_{k_2-1}&\ldots&c_{k_n-1}&{c_{m-1}}\\
\vdots&\vdots&\ &\vdots&\vdots\\
c_{k_1-n}&c_{k_2-n}&\ldots&c_{k_n-n}&{c_{m-n}}\\

\end{array}
\right|
$$
for $n=1,2,\ldots,r-1; k_1<k_2<\ldots<k_n<m,$ is nonnegative by hypothesis.
Since $\lim_{m\to \infty}c_m/c_{m-1}=1/T$, we have that
$\lim_{m\to \infty}c_m/c_{m-l}=1/T^l$.
Letting $m$ tend to $ \infty$ in the considered minors, we obtain

$$
\left|
\begin{array}{ccccc}
c_{k_1}&c_{k_2}&\ldots&c_{k_n}&1/T^n\\
c_{k_1-1}&c_{k_2-1}&\ldots&c_{k_n-1}&1/T^{n-1}\\
\vdots&\vdots&\ &\vdots&\vdots\\
c_{k_1-n}&c_{k_2-n}&\ldots&c_{k_n-n}&1\\

\end{array}
\right|\ge 0.
$$

Now, starting from the second row, we divide each row by
$T$ and subtract it from the previous, altering each row but the last.
We obtain

$$
0\le
\left|
\begin{array}{ccccc}
c^{(1)}_{k_1}&c^{(1)}_{k_2}&\ldots&c^{(1)}_{k_n}&0\\
c^{(1)}_{k_1-1}&c^{(1)}_{k_2-1}&\ldots&c^{(1)}_{k_n-1}&0\\
\vdots&\vdots&\ &\vdots&\vdots\\
c^{(1)}_{k_1-n+1}&c^{(1)}_{k_2-n+1}&\ldots&c^{(1)}_{k_n-n+1}&0\\
c^{(1)}_{k_1-n}&c^{(1)}_{k_2-n}&\ldots&c^{(1)}_{k_n-n}&1\\

\end{array}
\right|
$$
$$
=
\left|
\begin{array}{cccc}
c^{(1)}_{k_1}&c^{(1)}_{k_2}&\ldots&c^{(1)}_{k_n}\\
c^{(1)}_{k_1-1}&c^{(1)}_{k_2-1}&\ldots&c^{(1)}_{k_n-1}\\
\vdots&\vdots&\ &\vdots\\
c^{(1)}_{k_1-n+1}&c^{(1)}_{k_2-n+1}&\ldots&c^{(1)}_{k_n-n+1}\\

\end{array}
\right|.
$$

For proving that $\{c^{(1)}_k\}\in PF_{r-1}$ we introduce the
function $e^{\varepsilon z}=\sum_{k=0}^\infty b_kz^k$, and consider
$b_k=0$ for $k<0$. Note that (see~\cite{tp}, page 428)
$$
\left|
\begin{array}{cccc}
b_{k}&b_{k+1}&\ldots&b_{k+n}\\
b_{k-1}&b_{k}&\ldots&b_{k+n-1}\\
\vdots&\vdots&\ &\vdots\\
b_{k-n}&b_{k-n+1}&\ldots&b_{k}
\end{array}
\right|
>0
$$
for $k\ge 0$ and $n=0,1,2,\ldots.$

Let $g(z)=f_1(z)e^{\varepsilon z}=\sum_{k=0}^\infty a_kz^k$, where
$a_k=\sum_{m=0}^k c_m^{(1)}b_{k-m}=
\sum_{m=0}^k c_m^{(1)}\varepsilon^{k-m}/(k-m)!$, and let
$a_k=0$ for $k<0$.

By the Cauchy-Binet formula, we can write

$$
\left|
\begin{array}{cccc}
a_{k}&a_{k+1}&\ldots&a_{k+n}\\
a_{k-1}&a_{k}&\ldots&a_{k+n-1}\\
\vdots&\vdots&\ &\vdots\\
a_{k-n+1}&a_{k-n+2}&\ldots&a_{k}
\end{array}
\right|=
$$
$$
\sum_{m_1<\ldots<m_n}
\left|
\begin{array}{cccc}
c^{(1)}_{m_1}&c^{(1)}_{m_2}&\ldots&c^{(1)}_{m_n}\\
c^{(1)}_{m_1-1}&c^{(1)}_{m_2-1}&\ldots&c^{(1)}_{m_n-1}\\
\vdots&\vdots&\ &\vdots\\
c^{(1)}_{m_1-n+1}&c^{(1)}_{m_2-n+1}&\ldots&c^{(1)}_{m_n-n+1}

\end{array}
\right|
\left|
\begin{array}{cccc}
b_{k-m_1}&b_{k+1-m_1}&\ldots&b_{k+n-1-m_1}\\
b_{k-m_2}&b_{k+1-m_2}&\ldots&b_{k+n-1-m_2}\\
\vdots&\vdots&\ &\vdots\\
b_{k-m_n}&b_{k+1-m_n}&\ldots&b_{k+n-1-m_n}
\end{array}
\right|
$$

In the expression above each summand is nonnegative,
and the summand corresponding to
$m_1=0,m_2=1,\ldots,m_n=n-1$ is

$$
\left[c_0^{(1)}\right]^n
\left|
\begin{array}{cccc}
b_{k}&b_{k+1}&\ldots&b_{k+n-1}\\
b_{k-1}&b_{k}&\ldots&b_{k+n-2}\\
\vdots&\vdots&\ &\vdots\\
b_{k-n+1}&b_{k-n+2}&\ldots&b_{k}
\end{array}
\right|
>0
$$
for $k\ge 0$ and $n=0,1,2,\ldots$.

The series $\sum_{k=0}^\infty a_kz^k$ converges
for certain $\lambda,0<\lambda<1,$ and
the minors $\det\|a_{k+j-i}\|_{i,j=\overline{1,n}}, k\ge0,1\le n\le r-1$,
are strictly positive.
By Lemma 2, $\{a_k\}$ is a $PF_{r-1}$ sequence.
Also, $a_k\to c_k^{(1)}$, when $\varepsilon\to 0$.
Thus, $\{c_k^{(1)}\}\in PF_{r-1}$ as well.$\Box$

{\bf Proof of Corollary 1:}
By Theorem 4, the function $f_1(z)=(1-z/T)f(z)$
is a $PF_1$ g.f. and its radius of convergence is not less than $T$.
Thus, the limit in question
is either positive
or $+\infty$.$\Box$

{\bf Proof of Corollary 2:}
Applying $r-1$ times Theorem 4 (each
time the obtained function has radius of convergence equal to $T$),
we have that $(1-z/T)^{r-1}f(z)\in PF_1$
and the limit in question
is either positive
or $+\infty$.$\Box$

{\bf Acknowledgments:}

Deep gratitude to Prof. I.V.Ostrovskii and Prof. A.M.Vishnyakova
for fruitful discussions on the topic of this communication.

\end{document}